\documentclass[12pt,reqno]{amsart}

\usepackage[left=3.5cm,right=3.5cm]{geometry}
\usepackage{amssymb}
\usepackage{graphicx}
\usepackage{amscd}
\usepackage[pagebackref]{hyperref}
\usepackage{color}
\usepackage{tabularx}
\usepackage[table,dvipsnames]{xcolor}
\usepackage{float}
\usepackage{graphics,amsmath,amssymb}
\usepackage{amsthm}
\usepackage{amsfonts}
\usepackage{latexsym}
\usepackage{epsf}
\usepackage{xifthen}
\usepackage{mathrsfs}
\usepackage{dsfont}
\usepackage{makecell}
\usepackage{subfig}
\usepackage{amsmath}
\usepackage{listings}
\usepackage{etoolbox}
\usepackage{fancyhdr}
\usepackage{pdflscape}
\usepackage[title,toc,titletoc]{appendix}
\usepackage{enumitem}
\usepackage[noadjust]{cite}
\usepackage{tikz}
\usetikzlibrary{automata,positioning,arrows}
\usepackage{young}
\usepackage[object=vectorian]{pgfornament} 
\usepackage{lipsum,tikz}
\usepackage{multirow}
\usepackage[OT2,T1]{fontenc}
\usepackage{mathtools}

\hypersetup{
	colorlinks=true, 
	linktoc=all, 
	linkcolor=blue} 

\numberwithin{equation}{section}

\theoremstyle{plain}
\newtheorem{theorem}{Theorem}[section]
\newtheorem*{theorem*}{Theorem}

\newtheorem{corollary}[theorem]{Corollary}
\newtheorem{lemma}[theorem]{Lemma}
\newtheorem{proposition}[theorem]{Proposition}

\newtheorem{innercustomgeneric}{\customgenericname}
\providecommand{\customgenericname}{}
\newcommand{\newcustomtheorem}[2]{%
	\newenvironment{#1}[1]
	{%
		\renewcommand\customgenericname{#2}%
		\renewcommand\theinnercustomgeneric{##1}%
		\innercustomgeneric
	}
	{\endinnercustomgeneric}
}
\newcustomtheorem{ctheorem}{Theorem}
\newcustomtheorem{clemma}{Lemma}

\theoremstyle{definition}
\newtheorem{definition}[theorem]{Definition}

\newtheorem{example}[theorem]{Example}
\newtheorem*{example*}{Example}
\newtheorem*{examples*}{Examples}
\newtheorem{remark}[theorem]{Remark}
\newtheorem*{remark*}{Remark}
\newtheorem*{remarks*}{Remarks}
\newtheorem*{note*}{Note}

\newtheoremstyle{named}{}{}{\itshape}{}{\bfseries}{.}{.5em}{#1\thmnote{ #3}}
\theoremstyle{named}

\newtheoremstyle{customized}{}{}{\itshape}{}{\bfseries}{.}{.5em}{\thmnote{#3}}
\theoremstyle{customized}

\makeatletter
\patchcmd{\subsection}{\bfseries}{\bfseries\boldmath}{}{}
\makeatother

\DeclareMathAlphabet{\mydutchcal}{U}{dutchcal}{m}{n}

\newcommand{\cP}{\mydutchcal{P}}
\newcommand{\cD}{\mydutchcal{D}}

\newcommand{\cT}{\mydutchcal{T}}
\newcommand{\cU}{\mydutchcal{U}}
\newcommand{\cV}{\mydutchcal{V}}

\newcommand{\cPD}{\mydutchcal{PD}}
\newcommand{\cPDO}{\mydutchcal{PDO}}

\newcommand{\dis}{\mathsf{div}}
\newcommand{\len}{\mathsf{len}}

\newcommand{\even}{\mathsf{even}}

\newcommand{\UU}{\mathbf{U}}

\newcommand{\PDO}{\operatorname{PDO}}

\newcommand{\tmu}{\tilde{\mu}}
\newcommand{\tnu}{\tilde{\nu}}

\newcommand{\ttt}{\tilde{t}}

\newcommand{\doi}[1]{\href{https://dx.doi.org/#1}{DOI: #1}}

\newcommand{\OEIS}[1]{\href{https://oeis.org/#1}{#1}}

\title[Convolutive sequences,~III]{Convolutive sequences,~III:~A bijective approach to the $2$-convolutivity of the PDO function and beyond}

\author[S. Chern]{Shane Chern}
\address[S. Chern]{Fakult\"at f\"ur Mathematik, Universit\"at Wien, Oskar-Morgenstern-Platz 1, Wien 1090, Austria}
\email{chenxiaohang92@gmail.com, xiaohangc92@univie.ac.at}

\author[D. Eichhorn]{Dennis Eichhorn}
\address[D. Eichhorn]{Department of Mathematics, University of California Irvine, Irvine, CA 92697, USA}
\email{deichhor@math.uci.edu}

\author[S. Fu]{Shishuo Fu}
\address[S. Fu]{College of Mathematics and Statistics \& Center for Discrete Mathematics, Chongqing University, Chongqing 401331, China}
\email{fsshuo@cqu.edu.cn}

\author[J. A. Sellers]{James A. Sellers}
\address[J. A. Sellers]{Department of Mathematics, University of Minnesota Duluth, Duluth, MN 55811, USA}
\email{jsellers@d.umn.edu}

\date{}



\begin{document}
	
\sloppy

\begin{abstract}
	In 2002, Andrews, Lewis, and Lovejoy began the study of partitions with designated summands. A surprising result induced from their work is that the function $\operatorname{PDO}(n)$, which counts the number of such partitions of $n$ in which all parts are odd, satisfies the $2$-convolutive relation 
	$$\sum_{n\ge 0} \operatorname{PDO}(2n) q^n = \left(\sum_{n\ge 0} \operatorname{PDO}(n) q^n\right)^2.$$
	In previous work, the authors identified seven $2$-convolutive sequences, including $(\operatorname{PDO}(n))_{n\ge 0}$, such that each of them is generated by an eta-product. Moreover, the $2$-convolutivity of two such eta-products was proven combinatorially. However, the problem of finding a bijective proof for the $\operatorname{PDO}$ case remained open. In this paper, we resolve this problem with the long-desired proof based on Kim's bijection for a theta constant identity. Using the same technique, we also prove the $2$-convolutivity of another three known eta-products bijectively. Finally, we discover an overlooked statistic-preserving property of Kim's bijection, from which we obtain a combinatorial demonstration of a $2$-convolutive parametrization established in the second paper of the present series. 
\end{abstract}

\maketitle

\section{Introduction}

A \emph{partition}~\cite{And1998} of a natural number $n$ is a weakly decreasing sequence of positive integers, known as \emph{parts} of this partition, such that their sum equals $n$. Given a partition $\lambda$, we denote its \emph{weight}, written as $|\lambda|$, to be the sum of its parts. In addition, the weight of a partition tuple is the overall weight of all its components. Meanwhile, we introduce the statistic $\len(\lambda)$ to count the number of parts in $\lambda$.

In 2002, Andrews, Lewis, and Lovejoy~\cite{ALL2002} considered a new type of partitions, called \emph{partitions with designated summands}, by marking exactly one part of each different part size. For example, $6'+3'+3+2+2'+2$ is a partition with designated summands in which the only part of size $6$, the first part of size $3$, and the second part of size $2$ are marked. Let $\PDO(n)$ denote the number of such partitions of $n$ with the additional requirement that all parts are odd. A surprising result induced from the work of Andrews, Lewis, and Lovejoy is the relation~\cite[eq.~(1.6) and Theorem~21]{ALL2002}:
\begin{align}\label{eq:PDO-2-conv}
	\sum_{n\ge 0} \operatorname{PDO}(2n) q^n = \left(\sum_{n\ge 0} \operatorname{PDO}(n) q^n\right)^2.
\end{align}
This identity arose analytically from an explicit generating function expression for $\PDO(n)$ proven in \cite[eq.~(1.6)]{ALL2002}:
\begin{align}\label{eq:PDO-gf}
	\sum_{n\ge 0} \PDO(n) q^n = \frac{f_4 f_6^2}{f_1 f_3 f_{12}},
\end{align}
where
\begin{align*}
	f_i := (q^i;q^i)_\infty,
\end{align*}
with the conventional \emph{$q$-Pochhammer symbols} defined by
\begin{align*}
	(z;q)_\infty &:= \prod_{k\ge 0} (1-zq^k),\\
	(z_1,z_2,\ldots,z_l;q)_\infty &:= (z_1;q)_\infty (z_2;q)_\infty \cdots (z_l;q)_\infty.
\end{align*}
The expression on the right-hand side of \eqref{eq:PDO-gf} is often called an \emph{eta-product} because $f_1$ differs from the Dedekind eta function only by a prefactor $q^{\frac{1}{24}}$.

A sequence $(a(n))_{n\ge 0}$ satisfying the same relation as \eqref{eq:PDO-2-conv}, namely,
\begin{align}\label{eq:2-conv-def}
	\sum_{n\ge 0} a(2n) q^n = \left(\sum_{n\ge 0} a(n) q^n\right)^2,
\end{align}
was called \emph{$2$-convolutive} by the authors in \cite{CEFSI}. In the first two papers of the present series~\cite{CEFSI,CEFSII}, the $2$-convolutivity of seven eta-products, as listed in Table~\ref{tab:2-conv}, was established. However, most of these $2$-convolutive phenomena were shown in a purely analytic way, with only two of them having alternative combinatorial demonstrations. In particular, the problem of finding a bijective proof of the $2$-convolutivity of $\PDO(n)$ remained open.

\begin{table}[ht]
	\renewcommand{\arraystretch}{2.4} 
	\caption{List of $2$-convolutive eta-products in this series of papers}\label{tab:2-conv}
	\begin{center}
		\begin{tabular}{c|cccc} 
			\hline
			Label in \cite{CEFSII} & Eta-product & OEIS Entry & Analytic proof & Bijective proof\\
			\hline
			(II.1) & $\dfrac{f_2^6}{f_1^4 f_4^2}$ & \OEIS{A007096} & {\footnotesize\makecell{\cite{CEFSI}, Thm.~3.7\\\cite{CEFSII}, Thm.~1.1}} & {\footnotesize\makecell{\cite{CEFSI}, Thm.~2.10 \\ \textbf{This work, Thm.~\ref{th:A007096-bij}}}}\\
			(II.2) & $\dfrac{f_2^3 f_3}{f_1^3 f_6}$ & \OEIS{A128128} & {\footnotesize\makecell{\cite{CEFSI}, Thm.~3.6\\\cite{CEFSII}, Thms.~1.1/1.2}} & {\footnotesize\bfseries\makecell{This work\\ Thm.~\ref{th:A128128-bij}}}\\
			(II.3) & $\dfrac{f_2 f_4^2}{f_1^2 f_8}$ & \OEIS{A103258} & {\footnotesize\makecell{\cite{CEFSI}, Thm.~3.8\\\cite{CEFSII}, Thm.~1.1}} & {\footnotesize\makecell{\cite{CEFSI}, Thm.~2.11 \\ \textbf{This work, Thm.~\ref{th:A103258-bij}}}}\\
			(II.4) & $\dfrac{f_4 f_6^2}{f_1 f_3 f_{12}}$ & \OEIS{A102186} & {\footnotesize\makecell{\cite{CEFSI}, Thm.~3.4\\\cite{CEFSII}, Thms.~1.1/1.2}} & {\footnotesize\bfseries\makecell{This work\\ Thms.~\ref{thm:PDOb}/\ref{th:PDO-bij-step2}}}\\
			(II.5) & $\dfrac{f_1f_2f_6}{f_3^3}$ & \OEIS{A398619} & \cite{CEFSII}, Thm.~1.2 &\\
			(II.6) & $\dfrac{f_2f_3f_6}{f_1^2 f_9}$ & \OEIS{A398620} & \cite{CEFSII}, Thm.~1.2 &\\
			(II.7) & $\dfrac{f_6 f_{10}}{f_1 f_{15}}$ & \OEIS{A094023} & \cite{CEFSI}, Thm.~3.5 &\\[10pt]
			\hline
		\end{tabular}
	\end{center}
\end{table}

In \cite{SunKim1}, Sun Kim constructed a bijection for the general theta constant identity:
\begin{align*}
	&\langle q^{a_1};q^{2m}\rangle_\infty \langle q^{a_2};q^{2m}\rangle_\infty \langle q^{a_3};q^{2m}\rangle_\infty \langle q^{a_1+a_2+a_3};q^{2m}\rangle_\infty\\
	&\quad + \langle-q^{a_1};q^{2m}\rangle_\infty \langle-q^{a_2};q^{2m}\rangle_\infty \langle-q^{a_3};q^{2m}\rangle_\infty \langle-q^{a_1+a_2+a_3};q^{2m}\rangle_\infty\\
	&\quad\quad = 2(-q^{2m},-q^{2m};q^{2m})_\infty \langle-q^{a_1+a_2};q^{2m}\rangle_\infty \langle-q^{a_1+a_3};q^{2m}\rangle_\infty \langle-q^{a_2+a_3};q^{2m}\rangle_\infty,
\end{align*}
where
\begin{align*}
	\langle z;q\rangle_\infty := (z,q/z;q)_\infty,
\end{align*}
and $a_1$, $a_2$, and $a_3$ are odd positive integers no larger than $m$, with $a_1 + a_2 + a_3 <2m$. Her construction was motivated by earlier work of Warnaar~\cite{War2005} in which an identity originally discovered by Farkas and Kra~\cite{FK2000,FK2001} was studied combinatorially.

Now our key observation is that, by rewriting the eta-product in \eqref{eq:PDO-gf},
\begin{align*}
	\frac{f_4 f_6^2}{f_1 f_3 f_{12}} = (-q^2,-q^2;q^2)_\infty \langle -q;q^{12}\rangle_\infty \langle -q^{3};q^{12}\rangle_\infty \langle -q^{3};q^{12}\rangle_\infty \langle -q^{7};q^{12}\rangle_\infty,
\end{align*}
there are factors fitting well into Kim's setting. This then leads us to complete the long-desired bijective proof of the $2$-convolutivity of $\PDO(n)$.

\begin{theorem}\label{thm:PDOb}
	Let $\cPDO$ denote the set of partitions with designated summands in which all parts are odd, and let $[\cPDO]_\even$ denote the subset of $\cPDO$ consisting of all members with even weight. There is a weight-halving bijection from $[\cPDO]_\even$ to $\cPDO^2$. Therefore, the numbers $\PDO(n)$ form a $2$-convolutive sequence.
\end{theorem}

Kim's bijection also brings about combinatorial treatments of three more $2$-convolutive eta-products, namely, the series (II.1)--(II.3) in Table~\ref{tab:2-conv}. Interestingly, the $2$-convolutivity of these three eta-products, along with the generating series of $\PDO(n)$, can be parametrized by the relation in \cite[Theorem~1.1]{CEFSII}, with $X\mapsto -X$,
\begin{align}\label{eq:para-2-conv}
	\UU_2\left(\frac{f_2^2}{f_1^2} (-qX,-q/X;q^2)_\infty\right) = \left(\frac{f_2^2}{f_1^2} (-qX,-q/X;q^2)_\infty\right)^2,
\end{align}
where
\begin{align*}
	\UU_2\left(\sum_{n\ge 0} a(n) q^n\right) := \sum_{n\ge 0} a(2n) q^n.
\end{align*}
A combinatorial form of this parametrization related to  
$\cPDO$ 
was originally found by two of the authors in \cite{FuSel}. However, both proofs in \cite{CEFSII} and \cite{FuSel} are purely analytic.

It is then intriguing to know if Kim's bijection also handles the parametrization \eqref{eq:para-2-conv}. Since there is an extra parameter $X$, a combinatorial demonstration should involve a certain statistic for the objects in question. Unfortunately, Kim's original statement in \cite[Theorem~3.6]{SunKim1} contains no such statistic, making it unlikely to be used in a direct manner. On the other hand, a close examination of Kim's construction reveals an overlooked statistic-preserving property. In view of this new piece of information, a bijective proof of the parametrization \eqref{eq:para-2-conv} becomes possible, as will be shown in Theorem~\ref{th:para-bij}.

\textbf{Outline of the paper.} In Section~\ref{sec:prelim}, after introducing necessary notation and definitions, we explain the basic idea of our approach, from which the demand for Kim's bijection arises naturally. In Section~\ref{sec:BijPDO}, we present the long-awaited bijective proof of the $2$-convolutivity of $\PDO(n)$ claimed in Theorem~\ref{thm:PDOb}. Toward this proof, besides Kim's bijection, we also perform a combinatorial study of a generating function identity established by Andrews, Lewis, and Lovejoy. In Section~\ref{sec:BijProofs}, we apply the same technique to the eta-products (II.1)--(II.3) in Table~\ref{tab:2-conv} to demonstrate their $2$-convolutivity. In Section~\ref{sec:para}, we first review Kim's original construction of her bijection and then derive an overlooked statistic-preserving property, thereby completing a bijective treatment of the $2$-convolutive parametrization \eqref{eq:para-2-conv}. In Section~\ref{sec:conclusion}, we propose some questions for future study. Finally, in Appendix~\ref{app:example}, we offer a concrete example in full detail to illustrate the bijection stated in Theorem~\ref{thm:PDOb}.

\section{Preliminaries}\label{sec:prelim}

\subsection{Notation and definitions}

To study the combinatorics related to our $2$-convolutive sequences, we set forth the notation of several important sets of partitions or partition tuples:
\begin{itemize}[itemindent=*, leftmargin=*,itemsep=3pt]
	\item $\cP$: The set of partitions with no restrictions on the parts;
	
	\item $\cD$: The set of partitions into distinct parts;
	
	\item $\cD_{m,a}$: The set of partitions into distinct parts congruent to $a$ modulo $m$;
	
	\item $\cD_{m,\pm a}$: The Cartesian product $\cD_{m,a}\times \cD_{m,m-a}$;\footnote{We point out that this set is \emph{not} necessarily isomorphic to the set of partitions into distinct parts congruent to $\pm a$ modulo $m$; for example, $\cD_{m,\pm 0}$ is $\cD_{m,0}^2$ instead of $\cD_{m,0}$.}
	
	\item $\cD_{m,\pm \{a_1, a_2, \dots, a_l\} }$: The Cartesian product $\cD_{m,\pm a_1} \times \cD_{m,\pm a_2} \times \dots  \times \cD_{m,\pm a_l}$ with $\{a_1, a_2, \dots, a_l\}$ a list of integers that may repeat.
\end{itemize}

We will see from our proofs in Sections~\ref{sec:BijPDO} and \ref{sec:BijProofs} that many of the bijective manipulations are simply accomplished by rearranging components of a certain partition tuple. In view of this, the following definition is beneficial for the sake of brevity.

\begin{definition}\label{def:cong}
	Given two sets $\cT$ and $\cT'$ of partition tuples, if they can be bijectively transformed between one another by merely rearranging their components, then we write $\cT \cong \cT'$.
\end{definition}

As an example, by definition we have
\begin{align*}
	\cD_{m,\pm a} = \cD_{m,a}\times \cD_{m,m-a}
\end{align*}
and
\begin{align*}
	\cD_{m,\pm (m-a)} = \cD_{m,m-a}\times \cD_{m,a},
\end{align*}
so that there is a natural bijection between $\cD_{m,\pm a}$ and $\cD_{m,\pm (m-a)}$ by interchanging the two components in the ordered pairs. In accordance with Definition~\ref{def:cong}, we may write
\begin{align*}
	\cD_{m,\pm a}\cong \cD_{m,\pm (m-a)}.
\end{align*}

Now assume that $\cU$ is an arbitrary set of partitions or partition tuples.

\begin{definition}
	Let $k$ be a positive integer. We denote by $k\cU$ the set of partitions or partition tuples where each part in $\cU$ is multiplied by $k$, and by $\cU^k$ the Cartesian product
	\begin{align*}
		\cU^k := \underbrace{\cU\times \cdots \times \cU}_{\text{$k$ components}}.
	\end{align*}
\end{definition}

\begin{definition}\label{def:even}
	We denote by $[\cU]_\even$ the subset of $\cU$ consisting of all members with even weight.
\end{definition}

We record a simple observation.

\begin{lemma}\label{le:even-trivial}
	Let $\cU$ be an arbitrary set of partitions or partition tuples, and let $\cV$ be a set of partitions or partition tuples such that all its members have even weight. Then $[\cU\times \cV]_\even$ is identical to $[\cU]_\even\times \cV$.
\end{lemma}

\begin{proof}
	Note that the weight of the component from $\cV$ is always even. Thus, to pick out members from $\cU\times \cV$ with even weight, we must choose those such that the component from $\cU$ has even weight, thereby yielding the claimed $[\cU]_\even\times \cV$.
\end{proof}

\subsection{Basic idea in our approach}\label{sec:idea}

The operation given in Definition~\ref{def:even} plays a central role in our work. Suppose we want to show the $2$-convolutivity of a sequence $(u(n))_{n\ge 0}$, where $u(n)$ enumerates partition tuples of weight $n$ in a certain set $\cU$. That is,
\begin{align*}
	\sum_{n\ge 0} u(2n) q^n = \left(\sum_{n\ge 0}u(n) q^n\right)^2.
\end{align*}
Writing this relation combinatorially, we have
\begin{align*}
	\sum_{\pi\in [\cU]_\even} q^{\frac{|\pi|}{2}} = \left(\sum_{\sigma\in \cU} q^{|\sigma|}\right)^2 = \sum_{(\sigma,\rho) \in \cU^2} q^{|\sigma|+|\rho|}.
\end{align*}
Then it is sufficient to construct a weight-halving bijection from $[\cU]_\even$ to $\cU^2$.

For the moment, we concentrate on the situation where $\cU$ takes the form
\begin{align}\label{eq:U=U'*D}
	\cU' \times \cD_{2m,\pm\{2b_1,\ldots,2b_k\}}.
\end{align}
In light of Lemma~\ref{le:even-trivial},
\begin{align*}
	[\cU]_\even = [\cU']_\even \times \cD_{2m,\pm\{2b_1,\ldots,2b_k\}}.
\end{align*}
Ideally, we want $[\cU']_\even$ to possess a bijection with an explicitly expressed set of partition tuples. Since there is a weight-halving procedure following the bijection, we further hope that this set of partition tuples permits merely even parts; this discussion will be expanded in Section~\ref{sec:Kim-bij}. On the other hand, it is plain to see that $\cD_{2m,\pm\{2b_1,\ldots,2b_k\}}$ is in weight-halving bijection with $\cD_{m,\pm\{b_1,\ldots,b_k\}}$. However, the associated modulus of this set becomes $m$, while our desired $\cU^2$ satisfies
\begin{align*}
	\cU^2 \cong (\cU')^2 \times \cD_{2m,\pm\{2b_1,2b_1,\ldots,2b_k,2b_k\}},
\end{align*}
with the associated modulus of the second component in the above staying at $2m$. 

As such, if we want to keep track of the rearrangement of the components in a more transparent way, the following observation is useful in practice.

\begin{lemma}\label{le:trival}
	Let $m$ be a positive integer, and let $b_1,\ldots,b_k$ be integers. Then
	\begin{align}\label{eq:trivial-1}
		\cD_{2m,\pm\{2b_1,\ldots,2b_k\}}
	\end{align}
	is in weight-halving bijection with 
	\begin{align}\label{eq:trivial-2}
		\cD_{2m,\pm\{b_1,\ldots,b_k,m-b_k,\ldots,m-b_1\}}.
	\end{align}
\end{lemma}

\begin{proof}
	Note that all parts associated with the set \eqref{eq:trivial-1} are even. Dividing them by $2$ implies a weight-halving bijection from \eqref{eq:trivial-1} to $\cD_{m,\pm\{b_1,\ldots,b_k\}}$. We separate the parts in each component of the above according to their residues modulo $2m$. More precisely, for each $b$, we use the trivial bijection
	\begin{align*}
		\cD_{m,\pm b} \to \cD_{2m,\pm b} \times \cD_{2m,\pm (m-b)}.
	\end{align*}
	The resulting set is exactly \eqref{eq:trivial-2}.
\end{proof}

\subsection{Kim's bijection}\label{sec:Kim-bij}

As we have seen in Section~\ref{sec:idea}, there is a demand for a suitable set $\cU'$ in \eqref{eq:U=U'*D} such that it is possible to bijectively map $[\cU']_\even$ to a set of partition tuples in which only even parts are allowed. The cornerstone here is the aforementioned bijection of Kim, which can be rephrased using our notation in the following way.

\begin{theorem}[Kim, \cite{SunKim1}, Theorem 3.6]\label{thm:SunKim}
	Let $m$ be a positive integer, and let $a_1$, $a_2$, and $a_3$ be odd positive integers no larger than $m$, with $a_1 + a_2 + a_3 <2m$. 
	Then
	\begin{align}\label{eq:SK1}
		[\cD_{2m,\pm \{a_1,a_2,a_3,a_1+a_2+a_3\}}]_\even   
	\end{align}
	is in weight-preserving bijection with 
	\begin{align}\label{eq:SK2}
		\cD_{2m,\pm \{0,a_1+a_2,a_1+a_3,a_2+a_3\}}.
	\end{align}
\end{theorem}

\begin{remark}\label{rmk:Kim}
	Although we will use Kim's bijection in Sections~\ref{sec:BijPDO} and \ref{sec:BijProofs}, 
    we defer a detailed description of the construction to Section~\ref{sec:para}, from which we will see an overlooked statistic-preserving property. In addition, a concrete example will be illustrated in Appendix~\ref{ex:Kim}.
\end{remark}

Four of the seven known $2$-convolutive eta-products satisfy the condition that the $\cU'$ component as in \eqref{eq:U=U'*D} takes the form $\cD_{2m,\pm \{a_1,a_2,a_3,a_1+a_2+a_3\}}$ under the assumptions in Theorem~\ref{thm:SunKim}, and thus Kim's bijection allows us to bridge the gap between the relevant objects of the form \eqref{eq:SK1} and objects of the form \eqref{eq:SK2} in our bijective constructions. However, after performing the weight-halving action on \eqref{eq:SK2}, the associated modulus again reduces to $m$. To facilitate our analysis, the following corollary of Kim's bijection, which shares the same flavor as Lemma~\ref{le:trival}, is necessary.

\begin{corollary}\label{coro:Kim}
	Let $m$ be a positive integer, and let $a_1$, $a_2$, and $a_3$ be odd positive integers no larger than $m$, with $a_1 + a_2 + a_3 <2m$. Put
	\begin{align*}
		a_4 = \begin{cases}
			a_1+a_2+a_3, & \text{if $a_1+a_2+a_3< m$},\\
			2m-(a_1+a_2+a_3), & \text{if $a_1+a_2+a_3\ge m$}.
		\end{cases}
	\end{align*}
	Then
	\begin{align}\label{eq:Kim-coro-1}
		[\cD_{2m,\pm \{a_1,a_2,a_3,a_4\}}]_\even
	\end{align}
	is in weight-halving bijection with 
	\begin{align}\label{eq:Kim-coro-2}
		\cD_{2m,\pm\{0,\frac{a_1+a_2}{2},\frac{a_1+a_3}{2},\frac{a_2+a_3}{2},m-\frac{a_2+a_3}{2},m-\frac{a_1+a_3}{2},m-\frac{a_1+a_2}{2},m\}}.
	\end{align}
\end{corollary}

\begin{proof}
	Note that $\cD_{2m,\pm (2m-a_1-a_2-a_3)}\cong \cD_{2m,\pm (a_1+a_2+a_3)}$. Then \eqref{eq:Kim-coro-1} is in bijection with $[\cD_{2m,\pm \{a_1,a_2,a_3,a_1+a_2+a_3\}}]_\even$, which, by Kim's bijection, is further in weight-halving bijection with $\cD_{m,\pm\{0,\frac{a_1+a_2}{2},\frac{a_1+a_3}{2},\frac{a_2+a_3}{2}\}}$. Now separating the parts in each component of the above according to their residues modulo $2m$, we arrive at the desired set \eqref{eq:Kim-coro-2}.
\end{proof}

\section{The PDO function}\label{sec:BijPDO}

The objective of this section is the long-awaited bijective proof of the $2$-convolutivity of $\PDO(n)$, as stated in Theorem~\ref{thm:PDOb}, and we will fulfill it in two steps. Note that to illustrate a concrete example of our bijection, we often need to rearrange components in certain partition tuples, and to keep track of such manipulations, we have to color these components differently. In order to maintain the readability of the main text, we defer the example to Appendix~\ref{app:example}.

Let us start with a general generating function identity discovered by Andrews, Lewis, and Lovejoy~\cite{ALL2002}. For $S$ an arbitrary set of positive integers, let $\cPD_S$ be the set of partitions with designated summands such that all parts are from $S$. It was proven in \cite[Theorem~2.1]{ALL2002} that
\begin{align}\label{eq:ALL-gf}
	\sum_{\lambda\in \cPD_S} q^{|\lambda|} = \prod_{s\in S} \frac{1+q^{3s}}{(1-q^s)(1-q^{2s})}.
\end{align}

Now let $\cP_S$ and $\cD_S$ denote the subsets of $\cP$ and $\cD$, respectively, such that all parts are from $S$. We provide a bijective proof of \eqref{eq:ALL-gf}.

\begin{theorem}\label{th:ALL-bij}
	Let $S$ be an arbitrary set of positive integers. Then there is a weight-preserving bijection from $\cPD_S$ to $\cP_S\times 2\cP_S\times 3\cD_S$.
\end{theorem}

\begin{proof}
	Assume $\lambda\in \cPD_S$ is such that the part of size $s\in S$ appears $m_s$ times with the designation located at position $p_s$. Our purpose is to construct the corresponding partition triple $(\mu,\nu,\omega)\in \cP_S\times 2\cP_S\times 3\cD_S$. The basic principle is that for each $s\in S$, we map the $m_s$ parts of size $s$ to parts of sizes $s$, $2s$, and $3s$ in $\mu$, $\nu$, and $\omega$, respectively, such that the overall weight of these parts remains $m_s s$. If $m_s = 0$, then in $\mu$, $\nu$, and $\omega$, there are no parts of sizes $s$, $2s$, and $3s$, respectively. Now assume $m_s \ge 1$. The key observation is that by writing $m_s' := \lfloor\frac{m_s}{2}\rfloor + 1$, it is always true that $1\le m_s'\le m_s$. In addition,
	\begin{align*}
		m_s - m_s' = \left\lfloor\tfrac{m_s - 3}{2}\right\rfloor + 1.
	\end{align*}
	For the moment, there are two cases:
	\begin{itemize}[itemindent=*, leftmargin=*,itemsep=2pt]
		\item If $1\le p_s\le m_s'$, then we assign $m_s - 2(p_s-1)$ parts of size $s$ to $\mu$, $p_s-1$ parts of size $2s$ to $\nu$, and no part of size $3s$ to $\omega$. In this case, we have $0\le p_s-1 \le \lfloor\frac{m_s}{2}\rfloor$ so that $m_s - 2(p_s-1)$ is a nonnegative integer.
		
		\item If $m_s'+1 \le p_s\le m_s$, then we assign $m_s - 2(p_s-m_s'-1) - 3$ parts of size $s$ to $\mu$, $p_s-m_s'-1$ parts of size $2s$ to $\nu$, and one part of size $3s$ to $\omega$. In this case, if such $p_s$ exists, then $m_s\ge m_s'+1$ so that $m_s\ge 3$. Now, $0\le p_s-m_s'-1\le m_s - \lfloor\frac{m_s}{2}\rfloor - 2 = \lfloor\frac{m_s - 3}{2}\rfloor$, implying that $m_s - 2(p_s-m_s'-1) - 3$ is a nonnegative integer.
	\end{itemize}
	Therefore, we have an injective map from $\cPD_S$ to $\cP_S\times 2\cP_S\times 3\cD_S$.
	
	To see that the above map is surjective, we construct its inverse. For the triple $(\mu,\nu,\omega)\in \cP_S\times 2\cP_S\times 3\cD_S$, we assume that there are $m_{\mu,s}$ parts of size $s$ in $\mu$, $m_{\nu,s}$ parts of size $2s$ in $\nu$, and $m_{\omega,s}$ parts of size $3s$ in $\omega$, respectively, such that $m_{\mu,s}+2m_{\nu,s}+3m_{\omega,s} = m_s$. If $m_s = 0$, we have no part of size $s$ in $\lambda\in \cPD_S$. Now assuming $m_s\ge 1$, we assign $m_s$ parts of size $s$ to $\lambda\in \cPD_S$, and once again, there are two cases when designating these parts:
	\begin{itemize}[itemindent=*, leftmargin=*,itemsep=2pt]
		\item If $m_{\omega,s} = 0$, then $0\le m_{\nu,s}\le \lfloor\frac{m_s}{2}\rfloor$. Now the designated part of size $s$ in $\lambda$ is located at the position $p_s = m_{\nu,s} + 1$ so that $1\le p_s\le m_s'$.
		
		\item If $m_{\omega,s} = 1$, then $m_s\ge 3$ and $0\le m_{\nu,s}\le \lfloor\frac{m_s-3}{2}\rfloor$. Now the designated part of size $s$ in $\lambda$ is located at the position $p_s = m_s' + m_{\nu,s} + 1$ so that $m_s'+1 \le p_s\le m_s$.
	\end{itemize}
	The above construction produces the inverse of the previous map from $\cPD_S$ to $\cP_S\times 2\cP_S\times 3\cD_S$, thereby justifying its bijectivity.
\end{proof}

We may specialize Theorem~\ref{th:ALL-bij} to our desired set $\cPDO$ to get the first part of the bijective construction toward Theorem~\ref{thm:PDOb}.

\begin{corollary}\label{coro:PDO-bij-step1}
	There is a weight-preserving bijection from $\cPDO$ to $\cD\times \cD_{2,0}\times \cD_{6,3}$.
\end{corollary}

\begin{proof}
	Let $\mathbb{O}$ be the set of odd positive integers, and put $S = \mathbb{O}$ in Theorem~\ref{th:ALL-bij}. Then $\cPD_{\mathbb{O}}$ is indeed $\cPDO$, and it is in bijection with $\cP_{\mathbb{O}}\times 2\cP_{\mathbb{O}}\times 3\cD_{\mathbb{O}}$. In view of Glaisher's bijection~\cite{Gla1883} (see also \cite[Section~3.2.1]{Pak2006}),\footnote{Given an arbitrary partition in $\cP_{\mathbb{O}}$ with the odd part $i$ having multiplicity $m_i$, Glaisher's bijection maps it to a partition in $\cD$ such that it contains the part $i\cdot 2^r$ if and only if the $r$-th position of the binary expansion of $m_i$ has value $1$.} $\cP_{\mathbb{O}}$ is in bijection with $\cD$, and hence $2\cP_{\mathbb{O}}$ is in bijection with $2\cD$, which is exactly $\cD_{2,0}$. Finally, $\cD_{\mathbb{O}}$ is $\cD_{2,1}$ so that $3\cD_{\mathbb{O}}$ becomes $\cD_{6,3}$. The claimed bijection then follows.
\end{proof}

Our second ingredient for Theorem~\ref{thm:PDOb} is as follows.

\begin{theorem}\label{th:PDO-bij-step2}
	There is a weight-halving bijection from
	\begin{align*}
		[\cD\times \cD_{2,0}\times \cD_{6,3}]_\even
	\end{align*}
	to
	\begin{align*}
		(\cD\times \cD_{2,0}\times \cD_{6,3})^2.
	\end{align*}
\end{theorem}

\begin{remark}
	The statement of Theorem~\ref{th:PDO-bij-step2} is already sufficient to cast a bijective demonstration of the $2$-convolutivity of the sequence A102186 in OEIS~\cite{OEIS}, which corresponds to the eta-product (II.4) in Table~\ref{tab:2-conv}. This is because we can rewrite the eta-product as
	\begin{align*}
		\dfrac{f_4 f_6^2}{f_1 f_3 f_{12}} = (-q;q)_\infty (-q^2;q^2)_\infty (-q^3;q^6)_\infty,
	\end{align*}
	which is exactly the generating function for $\cD\times \cD_{2,0}\times \cD_{6,3}$.
\end{remark}

\begin{proof}
	We start by separating the parts in each component of the partition triples in $\cD\times \cD_{2,0}\times \cD_{6,3}$ according to their residues modulo $12$ and then rearranging the components. In particular, $\cD$ is in bijection with
	\begin{align*}
		\cD_{12,\pm\{1,2,3,4,5\}}\times \cD_{12,0}\times \cD_{12,6};
	\end{align*}
	$\cD_{2,0}$ is in bijection with
	\begin{align*}
		\cD_{12,\pm\{2,4\}}\times \cD_{12,0}\times \cD_{12,6};
	\end{align*}
	and $\cD_{6,3}$ is in bijection with $\cD_{12,\pm 3}$. We then combine the component $\cD_{12,0}$ from $\cD$ and the component $\cD_{12,0}$ from $\cD_{2,0}$ as the pair $\cD_{12,\pm 0}$. Similarly, joining $\cD_{12,6}$ from $\cD$ and $\cD_{12,6}$ from $\cD_{2,0}$ gives the pair $\cD_{12,\pm 6}$. Therefore, $\cD\times \cD_{2,0}\times \cD_{6,3}$ is in bijection with
	\begin{align*}
		\cD_{12,\pm \{0,1,2,2,3,3,4,4,5,6\}}\cong \cD_{12,\pm \{1,3,3,5\}}\times \cD_{12,\pm \{0,2,2,4,4,6\}}.
	\end{align*}
	Now in light of Lemma~\ref{le:even-trivial},
	\begin{align*}
		[\cD_{12,\pm \{1,3,3,5\}}\times \cD_{12,\pm \{0,2,2,4,4,6\}}]_\even = [\cD_{12,\pm \{1,3,3,5\}}]_\even \times \cD_{12,\pm \{0,2,2,4,4,6\}}.
	\end{align*}
	For the first set of partition tuples on the right, we use Corollary~\ref{coro:Kim} with $m=6$ and $(a_1,a_2,a_3)=(1,3,3)$ to get a weight-halving bijection
	\begin{align*}
		[\cD_{12,\pm \{1,3,3,5\}}]_\even \to \cD_{12,\pm \{0,2,2,3,3,4,4,6\}},
	\end{align*}
	and for the second set of partition tuples, we use Lemma~\ref{le:trival} to get a weight-halving bijection
	\begin{align*}
		\cD_{12,\pm \{0,2,2,4,4,6\}} \to \cD_{12,\pm \{0,1,1,2,2,3,3,4,4,5,5,6\}}.
	\end{align*}
	Finally, note that
	\begin{align*}
		\cD_{12,\pm \{0,2,2,3,3,4,4,6\}}\times \cD_{12,\pm \{0,1,1,2,2,3,3,4,4,5,5,6\}} \cong (\cD_{12,\pm \{0,1,2,2,3,3,4,4,5,6\}})^2.
	\end{align*}
	Since $\cD_{12,\pm \{0,1,2,2,3,3,4,4,5,6\}}$ is in bijection with $\cD\times \cD_{2,0}\times \cD_{6,3}$ as we have argued earlier, the desired weight-halving bijection is obtained.
\end{proof}

Now we are in a position to finalize the proof of Theorem~\ref{thm:PDOb}.

\begin{proof}[Proof of Theorem~\ref{thm:PDOb}]
	In light of Corollary~\ref{coro:PDO-bij-step1}, $[\cPDO]_\even$ is in weight-preserving bijection with $[\cD\times \cD_{2,0}\times \cD_{6,3}]_\even$. For this set of partition triples, Theorem~\ref{th:PDO-bij-step2} further asserts that it is in weight-halving bijection with $(\cD\times \cD_{2,0}\times \cD_{6,3})^2$. Then for the two components in this ordered pair, we use Corollary~\ref{coro:PDO-bij-step1} again to pull them back to $\cPDO$, thereby finishing the proof.
\end{proof}

\section{Three more $2$-convolutive sequences}\label{sec:BijProofs}

In a similar way, we can give bijective proofs of the $2$-convolutivity of three more sequences from Table~\ref{tab:2-conv}.
The first sequence, A128128 in OEIS that corresponds to the eta-product (II.2), had no bijective proof of $2$-convolutivity until now. The other two OEIS sequences, A007096 and A103258, which correspond to the eta-products (II.1) and (II.3), respectively, already have bijective proofs of their $2$-convolutivity given by the authors in \cite{CEFSI}. We include new combinatorial proofs below to demonstrate the power and versatility of the technique introduced in this work.

As in \cite[eq.~(1.3)]{CEFSI}, we define
\begin{align}\label{eq:eta-product}
	\sum_{n\ge 0} P_{m_1^{\delta_1}\cdots m_{J}^{\delta_J}}(n) q^n := \prod_{j=1}^J (q^{m_j};q^{m_j})_\infty^{\delta_j}.
\end{align}

\subsection{Sequence A128128}

The OEIS sequence A128128 is generated by
\begin{align*}
	\sum_{n\ge 0} P_{1^{-3} 2^3 3^{1} 6^{-1}}(n) q^n = \frac{f_2^3 f_{3}}{f_1^3 f_{6}} = (-q;q)_\infty (-q;q)_\infty (-q,-q^2;q^3)_\infty.
\end{align*}
In particular, this is the generating function for $\cD\times \cD\times \cD_{3,\pm 1}$.

\begin{theorem}\label{th:A128128-bij}
	There is a weight-halving bijection from
	\begin{align*}
		[\cD\times \cD\times \cD_{3,\pm 1}]_\even
	\end{align*}
	to
	\begin{align*}
		(\cD\times \cD\times \cD_{3,\pm 1})^2.
	\end{align*}
	Therefore, the numbers $P_{1^{-3} 2^3 3^{1} 6^{-1}}(n)$ form a $2$-convolutive sequence.
\end{theorem}

\begin{proof}
	By separating the parts in each component of the triples in
	\begin{align*}
		\cD\times \cD\times \cD_{3,\pm 1}
	\end{align*}
	according to their residues modulo $6$ and then rearranging the components, we see that this partition set is in bijection with
	\begin{align*}
		\cD_{6,\pm \{0,1,1,1,2,2,2,3\}}\cong \cD_{6,\pm \{1,1,1,3\}}\times \cD_{6,\pm \{0,2,2,2\}}.
	\end{align*}
	Now in light of Lemma~\ref{le:even-trivial},
	\begin{align*}
		[\cD_{6,\pm \{1,1,1,3\}}\times \cD_{6,\pm \{0,2,2,2\}}]_\even = [\cD_{6,\pm \{1,1,1,3\}}]_\even \times \cD_{6,\pm \{0,2,2,2\}}.
	\end{align*}
	For the first set of partition tuples on the right, we use Corollary~\ref{coro:Kim} with $m=3$ and $(a_1,a_2,a_3)=(1,1,1)$ to get a weight-halving bijection
	\begin{align*}
		[\cD_{6,\pm \{1,1,1,3\}}]_\even \to \cD_{6,\pm \{0,1,1,1,2,2,2,3\}},
	\end{align*}
	and for the second set of partition tuples, we use Lemma~\ref{le:trival} to get a weight-halving bijection
	\begin{align*}
		\cD_{6,\pm \{0,2,2,2\}} \to \cD_{6,\pm \{0,1,1,1,2,2,2,3\}}.
	\end{align*}
	Finally, note that
	\begin{align*}
		\cD_{6,\pm \{0,1,1,1,2,2,2,3\}}\times \cD_{6,\pm \{0,1,1,1,2,2,2,3\}} = (\cD_{6,\pm \{0,1,1,1,2,2,2,3\}})^2.
	\end{align*}
	Since $\cD_{6,\pm \{0,1,1,1,2,2,2,3\}}$ is in bijection with $\cD\times \cD\times \cD_{3,\pm 1}$ as we have shown at the beginning of this proof, the desired weight-halving bijection is obtained.
\end{proof}

\subsection{Sequence A007096}

The OEIS sequence A007096 is generated by
\begin{align*}
	\sum_{n\ge 0} P_{1^{-4} 2^6 4^{-2}}(n) q^n = \frac{f_2^6}{f_1^4 f_4^2} = (-q;q)_\infty (-q;q)_\infty (-q,-q;q^2)_\infty.
\end{align*}
In particular, this is the generating function for $\cD\times \cD\times \cD_{2,\pm 1}$.

\begin{theorem}
	\label{th:A007096-bij}
	There is a weight-halving bijection from
	\begin{align*}
		[\cD\times \cD\times \cD_{2,\pm 1}]_\even
	\end{align*}
	to
	\begin{align*}
		(\cD\times \cD\times \cD_{2,\pm 1})^2.
	\end{align*}
	Therefore, the numbers $P_{1^{-4} 2^6 4^{-2}}(n)$ form a $2$-convolutive sequence.
\end{theorem}

\begin{remark}
	The set of partition tuples considered in this theorem differs from the one in \cite[Theorem~2.10]{CEFSI}, which deals with the same sequence $P_{1^{-4} 2^6 4^{-2}}(n)$.
\end{remark}

\begin{proof}
	By separating the parts in each component of the triples in
	\begin{align*}
		\cD\times \cD\times \cD_{2,\pm 1}
	\end{align*}
	according to their residues modulo $4$ and then rearranging the components, we see that this partition set is in bijection with
	\begin{align*}
		\cD_{4,\pm \{0,1,1,1,1,2\}}\cong \cD_{4,\pm \{1,1,1,1\}}\times \cD_{4,\pm \{0,2\}}.
	\end{align*}
	Now in light of Lemma~\ref{le:even-trivial},
	\begin{align*}
		[\cD_{4,\pm \{1,1,1,1\}}\times \cD_{4,\pm \{0,2\}}]_\even = [\cD_{4,\pm \{1,1,1,1\}}]_\even \times \cD_{4,\pm \{0,2\}}.
	\end{align*}
	For the first set of partition tuples on the right, we use Corollary~\ref{coro:Kim} with $m=2$ and $(a_1,a_2,a_3)=(1,1,1)$ to get a weight-halving bijection
	\begin{align*}
		[\cD_{4,\pm \{1,1,1,1\}}]_\even \to \cD_{4,\pm \{0,1,1,1,1,1,1,2\}},
	\end{align*}
	and for the second set of partition tuples, we use Lemma~\ref{le:trival} to get a weight-halving bijection
	\begin{align*}
		\cD_{4,\pm \{0,2\}} \to \cD_{4,\pm \{0,1,1,2\}}.
	\end{align*}
	Finally, note that
	\begin{align*}
		\cD_{4,\pm \{0,1,1,1,1,1,1,2\}}\times \cD_{4,\pm \{0,1,1,2\}} \cong (\cD_{4,\pm \{0,1,1,1,1,2\}})^2.
	\end{align*}
	Since $\cD_{4,\pm \{0,1,1,1,1,2\}}$ is in bijection with $\cD\times \cD\times \cD_{2,\pm 1}$ as we have shown at the beginning of this proof, the desired weight-halving bijection is obtained.
\end{proof}

\subsection{Sequence A103258}

The OEIS sequence A103258 is generated by
\begin{align*}
	\sum_{n\ge 0} P_{1^{-2} 2^1 4^{2} 8^{-1}}(n) q^n = \frac{f_2 f_4^2}{f_1^2 f_8} = (-q;q)_\infty (-q;q)_\infty (-q^2,-q^6;q^8)_\infty.
\end{align*}
In particular, this is the generating function for $\cD\times \cD\times \cD_{8,\pm 2}$.

\begin{theorem}
	\label{th:A103258-bij}
	There is a weight-halving bijection from
	\begin{align*}
		[\cD\times \cD\times \cD_{8,\pm 2}]_\even
	\end{align*}
	to
	\begin{align*}
		(\cD\times \cD\times \cD_{8,\pm 2})^2.
	\end{align*}
	Therefore, the numbers $P_{1^{-2} 2^1 4^{2} 8^{-1}}(n)$ form a $2$-convolutive sequence.
\end{theorem}

\begin{remark}
	The set of partition tuples considered in this theorem differs from the one in \cite[Theorem~2.11]{CEFSI}, which deals with the same sequence $P_{1^{-2} 2^1 4^{2} 8^{-1}}(n)$.
\end{remark}

\begin{proof}
	By separating the parts in each component of the triples in
	\begin{align*}
		\cD\times \cD\times \cD_{8,\pm 2}
	\end{align*}
	according to their residues modulo $8$ and then rearranging the components, we see that this partition set is in bijection with
	\begin{align*}
		\cD_{8,\pm \{0,1,1,2,2,2,3,3,4\}}\cong \cD_{8,\pm \{1,1,3,3\}}\times \cD_{8,\pm \{0,2,2,2,4\}}.
	\end{align*}
	Now in light of Lemma~\ref{le:even-trivial},
	\begin{align*}
		[\cD_{8,\pm \{1,1,3,3\}}\times \cD_{8,\pm \{0,2,2,2,4\}}]_\even = [\cD_{8,\pm \{1,1,3,3\}}]_\even \times \cD_{8,\pm \{0,2,2,2,4\}}.
	\end{align*}
	For the first set of partition tuples on the right, we use Corollary~\ref{coro:Kim} with $m=4$ and $(a_1,a_2,a_3)=(1,1,3)$ to get a weight-halving bijection
	\begin{align*}
		[\cD_{8,\pm \{1,1,3,3\}}]_\even \to \cD_{8,\pm \{0,1,2,2,2,2,3,4\}},
	\end{align*}
	and for the second set of partition tuples, we use Lemma~\ref{le:trival} to get a weight-halving bijection
	\begin{align*}
		\cD_{8,\pm \{0,2,2,2,4\}} \to \cD_{8,\pm \{0,1,1,1,2,2,3,3,3,4\}}.
	\end{align*}
	Finally, note that
	\begin{align*}
		\cD_{8,\pm \{0,1,2,2,2,2,3,4\}}\times \cD_{8,\pm \{0,1,1,1,2,2,3,3,3,4\}} = (\cD_{8,\pm \{0,1,1,2,2,2,3,3,4\}})^2.
	\end{align*}
	Since $\cD_{8,\pm \{0,1,1,2,2,2,3,3,4\}}$ is in bijection with $\cD\times \cD\times \cD_{8,\pm 2}$ as we have shown at the beginning of this proof, the desired weight-halving bijection is obtained.
\end{proof}

\section{The first Chern--Eichhorn--Fu--Sellers parametrization}\label{sec:para}

Two parametrizations for the $2$-convolutive eta-products were established by the authors in \cite{CEFSII}. In this work, we focus on the first one, which has been given in \eqref{eq:para-2-conv}. Reformulating the bivariate series in this $2$-convolutive relation, we have
\begin{align}
	\frac{f_2^2}{f_1^2} (-qX,-q/X;q^2)_\infty = (-q,-q;q)_\infty (-qX,-q/X;q^2)_\infty.
\end{align}

Now given $\cD_{m,\pm a} = \cD_{m,a}\times \cD_{m,m-a}$, we always use $\kappa_\pm$ to denote an ordered pair $(\kappa_+,\kappa_-)$ in this set. By abuse of notation, we write
\begin{align*}
	|\kappa_\pm| := |\kappa_+| + |\kappa_-|,
\end{align*}
and
\begin{align*}
	\len(\kappa_\pm) := \len(\kappa_+) - \len(\kappa_-).
\end{align*}
It is clear that for $1\le a\le m-1$,
\begin{align*}
	\sum_{\kappa_\pm\in \cD_{m,\pm a}} X^{\len(\kappa_\pm)} q^{|\kappa_\pm|} = (-q^a X,-q^{m-a}/X ;q^m)_\infty.
\end{align*}
Therefore,
\begin{align}
	\sum_{(\pi_\pm^{(1)}, \pi_\pm^{(2)}) \in \cD_{1,\pm 0} \times \cD_{2,\pm 1}} X^{\len(\pi_\pm^{(2)})} q^{|\pi_\pm^{(1)}| + |\pi_\pm^{(2)}|} = (-q,-q;q)_\infty (-qX,-q/X;q^2)_\infty.
\end{align}
As such, our goal is to prove the following combinatorial recasting of \eqref{eq:para-2-conv}.

\begin{theorem}\label{th:para-bij}
	There is a weight-halving bijection
	\begin{align*}
		\begin{array}{ccc}
			[\cD_{1,\pm 0} \times \cD_{2,\pm 1}]_\even & \to & (\cD_{1,\pm 0} \times \cD_{2,\pm 1})^2\\[3pt]
			(\pi_\pm^{(1)}, \pi_\pm^{(2)}) & \mapsto & \big((\sigma_\pm^{(1)}, \sigma_\pm^{(2)}), (\rho_\pm^{(1)}, \rho_\pm^{(2)})\big)
		\end{array}
	\end{align*}
	such that
	\begin{align}
		\len(\pi_\pm^{(2)}) = \len(\sigma_\pm^{(2)}) + \len(\rho_\pm^{(2)}).
	\end{align}
\end{theorem}

\subsection{Kim's bijection revisited}\label{sec:Kim-bij-2}

As promised in Remark~\ref{rmk:Kim}, we will give a detailed description of Kim's bijection. To do so, we rephrase her original construction as a quintipartite procedure. Throughout, assume that
\begin{align*}
	\big(\mu_\pm^{(1)}, \mu_\pm^{(2)}, \mu_\pm^{(3)}, \mu_\pm^{(4)}\big)
\end{align*}
is a quadruple in $[\cD_{2m,\pm \{a_1,a_2,a_3,a_1+a_2+a_3\}}]_\even$, using the notation of Theorem~\ref{thm:SunKim}.

\begin{enumerate}[label={\textup{(\arabic*)~}},leftmargin=*,labelsep=0cm,align=left]
	\item Standardize our quadruple modulo $2m$ to
	\begin{align*}
		\big(\tmu_\pm^{(1)}, \tmu_\pm^{(2)}, \tmu_\pm^{(3)}, \tmu_\pm^{(4)}\big).
	\end{align*}
	Here we put
	\begin{align*}
		\tmu_+^{(i)} := \left\lfloor\tfrac{\mu_+^{(i)}}{2m}\right\rfloor
	\end{align*}
	in the sense that we divide each part in $\mu_+^{(i)}$ by $2m$ and then take the floor of this value --- in particular, we allow $0$ as a part in $\tmu_+^{(i)}$. In a similar manner,
	\begin{align*}
		\tmu_-^{(i)} := \left\lceil\tfrac{\mu_-^{(i)}}{2m}\right\rceil
	\end{align*}
	with the ceiling function utilized instead.
	
	\item Apply Wright's bijection~\cite{Wright} (see also \cite[Section~6.2.1]{Pak2006}) to map each $\tmu_\pm^{(i)}$ to an unrestricted partition $\lambda^{(i)}$ joint with an integer $t^{(i)}$ so that we have
	\begin{align*}
		\big((\lambda^{(1)}, \lambda^{(2)}, \lambda^{(3)}, \lambda^{(4)}),\, (t^{(1)}, t^{(2)}, t^{(3)}, t^{(4)})\big).
	\end{align*}
	Note that Wright's bijection explains the identity
	\begin{align*}
		(-1,-q;q)_\infty = \frac{1}{(q;q)_\infty} \sum_{t=-\infty}^{\infty} q^{\binom{t}{2}},
	\end{align*}
	so that in our notation, with the superscript ``$(i)$'' removed for simplicity, there holds
	\begin{align*}
		|\tmu_\pm| = |\tmu_+|+|\tmu_-| = |\lambda| + \tbinom{t}{2}.
	\end{align*}
	The value of $t$ is given by
	\begin{align*}
		t:= \len(\tmu_\pm).
	\end{align*}
	If $t\ge 1$, then we split $\tmu_+$ into an unrestricted partition and a staircase partition $(\len(\tmu_+)-1)+\cdots+1+0$. Next, for $\tmu_-$, we add its $i$-th part to $t-1+i$ for each $i$. Finally, to the end of this new partition, we append the conjugate of the unrestricted partition separated from $\tmu_+$ to derive the desired unrestricted partition $\lambda$. In the staircase partition from $\tmu_+$, the remaining parts are $(t-1)+\cdots+1+0$, the weight of which is the desired $\binom{t}{2}$. The following manipulations of Ferrers diagrams illustrate this process:
	\begin{align*}
		\vcenter{\hbox{\begin{tikzpicture}[scale=0.45]
					\fill[color=red] foreach \Z [count=\Y] in {8,7,6,4,3,1} {foreach \X in {1,...,\Z} {(\X/2,-\Y/2) circle[radius=3pt]}};
					\fill[color=blue,xshift=5cm] foreach \Z [count=\Y] in {5,4,1} {foreach \X in {1,...,\Z} {(\X/2,-\Y/2) circle[radius=3pt]}};
		\end{tikzpicture}}}
		\quad\mapsto\quad
		\vcenter{\hbox{\begin{tikzpicture}[scale=0.45]
					\fill[color=red] foreach \Z [count=\Y] in {5,4,3,2,1} {foreach \X in {1,...,\Z} {(-\X/2+1/2,-\Y/2) circle[radius=3pt]}};
					\fill[color=olive] foreach \Z [count=\Y] in {3,3,3,2,2,1} {foreach \X in {1,...,\Z} {(\X/2,-\Y/2) circle[radius=3pt]}};
					\fill[color=teal] foreach \Z [count=\Y] in {2,1} {foreach \X in {1,...,\Z} {(-\X/2-1,-\Y/2) circle[radius=3pt]}};
					\fill[color=blue,xshift=2.5cm] foreach \Z [count=\Y] in {5,4,1} {foreach \X in {1,...,\Z} {(\X/2,-\Y/2) circle[radius=3pt]}};
		\end{tikzpicture}}}
		\quad\mapsto\quad
		\vcenter{\hbox{\begin{tikzpicture}[scale=0.45]
					\fill[color=olive] foreach \Z [count=\Y] in {8,8,6,6,5,3} {foreach \X in {1,...,\Z} {(\X/2,-\Y/2) circle[radius=3pt]}};
					\fill[color=red] foreach \Z [count=\Y] in {8,8,6} {foreach \X in {1,...,\Z} {(\X/2,-\Y/2) circle[radius=3pt]}};
					\fill[color=blue] foreach \Z [count=\Y] in {5,4,1} {foreach \X in {1,...,\Z} {(\X/2,-\Y/2) circle[radius=3pt]}};
					\fill[color=teal,xshift=5cm] foreach \Z [count=\Y] in {2,1} {foreach \X in {1,...,\Z} {(\X/2,-\Y/2) circle[radius=3pt]}};
		\end{tikzpicture}}}\quad.
	\end{align*}
	If $t\le 0$, then we simply append a part of size $0$ to $\tmu_-$ and switch its role with $\tmu_+$.
	
	\item Map the quadruple $(t^{(1)}, t^{(2)}, t^{(3)}, t^{(4)})$ to $(\ttt^{(1)}, \ttt^{(2)}, \ttt^{(3)}, \ttt^{(4)})$ such that they satisfy the following parametrizations (see Remark~\ref{rmk:i is integer}) via $(i,j,k,l)\in\mathbb{Z}^4$:
	\begin{align*}
		\begin{array}{lp{0.5em}lp{0.5em}lp{0.5em}l}
			t^{(1)}=-2i-j+k+l,&& t^{(2)}=k-j,&& t^{(3)}=l-j,&& t^{(4)}=j;\\
			\ttt^{(1)}=-2j-i+k+l,&& \ttt^{(2)}=k-i,&& \ttt^{(3)}=l-i,&& \ttt^{(4)}=i.
		\end{array}
	\end{align*}
	Then we have
	\begin{align*}
		\big((\lambda^{(1)},\lambda^{(2)},\lambda^{(3)},\lambda^{(4)}),\, (\ttt^{(1)}, \ttt^{(2)}, \ttt^{(3)}, \ttt^{(4)})\big).
	\end{align*}
	
	\item Apply the inverse of Wright's map to each $\lambda^{(i)}$ and $\ttt^{(i)}$ to get an ordered pair $\tnu_\pm^{(i)}$, so that we have the quadruple
	\begin{align*}
		\big(\tnu_\pm^{(1)}, \tnu_\pm^{(2)}, \tnu_\pm^{(3)}, \tnu_\pm^{(4)}\big).
	\end{align*}
    This is the reverse of Step (2).

	\item  Map each $\tnu_\pm^{(i)}$ to the ordered pair $(2m\tnu_+^{(i)}+e_+^{(i)},2m\tnu_-^{(i)}+e_-^{(i)})$, where $2m\tnu+e$ means multiplying each part of $\tnu$ by $2m$, then adding $e$, and finally removing possible parts of size $0$. Here we put\footnote{We have corrected misprints on \cite[p.~709]{SunKim1} regarding the values of $e_\pm^{(3)}$ and $e_\pm^{(4)}$.}
	\begin{align*}
		(e_+^{(i)},e_-^{(i)}) := \begin{cases}
			(0,0), & \text{if $i=1$},\\
			(a_1+a_2,-a_1-a_2), & \text{if $i=2$},\\
			(a_1+a_3,-a_1-a_3), & \text{if $i=3$},\\
			(a_2+a_3,-a_2-a_3), & \text{if $i=4$}.
		\end{cases}
	\end{align*}
	This is the reverse of Step (1) but with a different set of residues used. Now we arrive at the desired quadruple
	\begin{align*}
		\big(\nu_\pm^{(1)}, \nu_\pm^{(2)}, \nu_\pm^{(3)}, \nu_\pm^{(4)}\big),
	\end{align*}
	which is in $\cD_{2m,\pm \{0,a_1+a_2,a_1+a_3,a_2+a_3\}}$.
\end{enumerate}

\begin{remark}\label{rmk:i is integer}
    In Step~(3) of the above procedure, one gets from the parametrization of the $t$-quadruple that $i=\frac{1}{2}(-t^{(1)}+t^{(2)}+t^{(3)}+t^{(4)})$, which in general is not necessarily an integer. In our case, however, since all parts in the partition quadruple $(\mu_{\pm}^{(1)},\mu_{\pm}^{(2)},\mu_{\pm}^{(3)},\mu_{\pm}^{(4)})$ from $[\cD_{2m,\pm \{a_1,a_2,a_3,a_1+a_2+a_3\}}]_\even$ are odd so that the number of parts has the same parity as the partition weight, we then have
    \begin{align*}
    	2i \equiv \sum_{i=1}^4 t^{(i)}= \sum_{i=1}^4\len(\tmu_{\pm}^{(i)}) = \sum_{i=1}^4\len(\mu_{\pm}^{(i)})\equiv\sum_{i=1}^4|\mu_{\pm}^{(i)}|\equiv 0\pmod 2.
    \end{align*}
    Thus, $i$ is indeed an integer, and the parametrizations via $(i,j,k,l)\in\mathbb{Z}^4$ are valid.
\end{remark}

Although Kim did not convey it in her statement of \cite[Theorem~3.6]{SunKim1}, it can be seen from Step~(3) of the above procedure that
\begin{align*}
	t^{(2)} + t^{(4)} &= \ttt^{(2)} + \ttt^{(4)},\\
	t^{(3)} + t^{(4)} &= \ttt^{(3)} + \ttt^{(4)}.
\end{align*}
In addition, Wright's map tells us that
\begin{align*}
	&t^{(2)} = \len(\tmu_\pm^{(2)}) = \len(\mu_\pm^{(2)}),\quad && \ttt^{(2)} = \len(\tnu_\pm^{(2)}) = \len(\nu_\pm^{(2)}),\\
	&t^{(3)} = \len(\tmu_\pm^{(3)}) = \len(\mu_\pm^{(3)}),\quad && \ttt^{(3)} = \len(\tnu_\pm^{(3)}) = \len(\nu_\pm^{(3)}),\\
	&t^{(4)} = \len(\tmu_\pm^{(4)}) = \len(\mu_\pm^{(4)}),\quad && \ttt^{(4)} = \len(\tnu_\pm^{(4)}) = \len(\nu_\pm^{(4)}).
\end{align*}
The identities in the right column are correct because for $i\in \{2,3,4\}$, no part of size $0$ occurs when producing $\nu_\pm^{(i)}$ from $\tnu_\pm^{(i)}$, in light of the fact that $e_\pm^{(i)}$ are nonzero. Our observation can be summarized as the following overlooked statistic-preserving property.

\begin{proposition}\label{prop:Kim-stat}
	Kim's bijection in Theorem~\ref{thm:SunKim},
	\begin{align*}
		\begin{array}{ccc}
			[\cD_{2m,\pm \{a_1,a_2,a_3,a_1+a_2+a_3\}}]_\even & \to & \cD_{2m,\pm \{0,a_1+a_2,a_1+a_3,a_2+a_3\}}\\[3pt]
			(\mu_\pm^{(1)}, \mu_\pm^{(2)}, \mu_\pm^{(3)}, \mu_\pm^{(4)}) & \mapsto & (\nu_\pm^{(1)}, \nu_\pm^{(2)}, \nu_\pm^{(3)}, \nu_\pm^{(4)})
		\end{array}
	\end{align*}
	satisfies
	\begin{align}
		\len(\mu_\pm^{(2)}) + \len(\mu_\pm^{(4)}) = \len(\nu_\pm^{(2)}) + \len(\nu_\pm^{(4)}),
	\end{align}
	and
	\begin{align}
		\len(\mu_\pm^{(3)}) + \len(\mu_\pm^{(4)}) = \len(\nu_\pm^{(3)}) + \len(\nu_\pm^{(4)}).
	\end{align}
\end{proposition}

\subsection{Proof of Theorem~\ref{th:para-bij}}

Now we prove Theorem~\ref{th:para-bij}. Recall that we have assumed
\begin{align*}
	\big(\pi_\pm^{(1)}, \pi_\pm^{(2)}\big)\in [\cD_{1,\pm 0} \times \cD_{2,\pm 1}]_\even.
\end{align*}
Since $\cD_{1,\pm 0}$ is in bijection with $\cD_{4,\pm 1}\times \cD_{4,\pm 1}\times \cD_{2,\pm 0}$, we map $\pi_\pm^{(1)}$ to the triple
\begin{align*}
	\big(\mu_\pm^{(1)}, \mu_\pm^{(2)}, \xi_\pm\big) \in \cD_{4,\pm 1}\times \cD_{4,\pm 1}\times \cD_{2,\pm 0}.
\end{align*}
Meanwhile, we know that $\cD_{2,1}$ is in bijection with $\cD_{4,1}\times \cD_{4,3}$. Thus, for $\cD_{2,\pm 1} = \cD_{2,1}\times \cD_{2,1}$, we map $\pi_+^{(2)}$ from the first $\cD_{2,1}$ to the ordered pair $(\kappa_+^{(3)}, \kappa_+^{(4)})$ and $\pi_-^{(2)}$ from the second $\cD_{2,1}$ to the ordered pair $(\kappa_-^{(3)}, \kappa_-^{(4)})$, both in $\cD_{4,1}\times \cD_{4,3}$. Rearranging them, we then get the bijection
\begin{align*}
	\begin{array}{ccc}
		\cD_{2,\pm 1} & \to & \cD_{4,\pm 1}\times \cD_{4,\pm 3}\\[3pt]
		\pi_\pm^{(2)} & \mapsto & (\mu_\pm^{(3)}, \mu_\pm^{(4)})
	\end{array}
\end{align*}
where we put
\begin{align*}
	\mu_\pm^{(3)} := \big(\kappa_+^{(3)}, \kappa_-^{(4)}\big),\qquad\qquad \mu_\pm^{(4)} := \big(\kappa_+^{(4)}, \kappa_-^{(3)}\big).
\end{align*}
Thus, we arrive at the bijection
\begin{align*}
	\begin{array}{ccc}
		\cD_{1,\pm 0} \times \cD_{2,\pm 1} & \to & \cD_{4,\pm \{1,1,1,3\}}\times \cD_{2,\pm 0}\\[3pt]
		(\pi_\pm^{(1)}, \pi_\pm^{(2)}) & \mapsto & \big((\mu_\pm^{(1)}, \mu_\pm^{(2)}, \mu_\pm^{(3)}, \mu_\pm^{(4)}), \xi_\pm\big)
	\end{array}
\end{align*}
such that
\begin{align*}
	\len(\mu_\pm^{(3)}) + \len(\mu_\pm^{(4)}) &= \big(\len(\kappa_+^{(3)}) - \len(\kappa_-^{(4)})\big) + \big(\len(\kappa_+^{(4)}) - \len(\kappa_-^{(3)})\big)\\
	&= \big(\len(\kappa_+^{(3)}) + \len(\kappa_+^{(4)})\big) - \big(\len(\kappa_-^{(3)}) + \len(\kappa_-^{(4)})\big)\\
	&= \len(\pi_+^{(2)}) - \len(\pi_-^{(2)})\\
	&= \len(\pi_\pm^{(2)}).
\end{align*}
In light of Lemma~\ref{le:even-trivial},
\begin{align*}
	[\cD_{4,\pm \{1,1,1,3\}}\times \cD_{2,\pm 0}]_\even = [\cD_{4,\pm \{1,1,1,3\}}]_\even \times \cD_{2,\pm 0}.
\end{align*}
Halving the second component $\cD_{2,\pm 0}$ gives
\begin{align*}
	\begin{array}{ccc}
		\cD_{2,\pm 0} & \xrightarrow{\text{halve}} & \cD_{1,\pm 0}\\[3pt]
		\xi_\pm & \xmapsto{\text{\phantom{halve}}} & \rho_\pm^{(1)}
	\end{array}.
\end{align*}
We then apply Theorem~\ref{thm:SunKim} and Proposition~\ref{prop:Kim-stat} to the first component $[\cD_{4,\pm \{1,1,1,3\}}]_\even$ and halve the resulting quadruple to get
\begin{align*}
	\begin{array}{ccc}
		[\cD_{4,\pm \{1,1,1,3\}}]_\even & \xrightarrow{\text{halve}} & \cD_{2,\pm \{0,1,1,1\}}\\[3pt]
		(\mu_\pm^{(1)}, \mu_\pm^{(2)}, \mu_\pm^{(3)}, \mu_\pm^{(4)}) & \xmapsto{\text{\phantom{halve}}} & (\omega_\pm^{(1)}, \omega_\pm^{(2)}, \sigma_\pm^{(2)}, \rho_\pm^{(2)})
	\end{array}
\end{align*}
such that
\begin{align*}
	\len(\sigma_\pm^{(2)}) + \len(\rho_\pm^{(2)}) = \len(\mu_\pm^{(3)}) + \len(\mu_\pm^{(4)}) = \len(\pi_\pm^{(2)}).
\end{align*}
Finally, we utilize the bijection
\begin{align*}
	\begin{array}{ccc}
		\cD_{2,\pm 0}\times \cD_{2,\pm 1} & \to & \cD_{1,\pm 0}\\
		(\omega_\pm^{(1)}, \omega_\pm^{(2)}) & \mapsto & \sigma_\pm^{(1)}
	\end{array}.
\end{align*}
Putting all procedures together, we have
\begin{equation*}
	\resizebox{\textwidth}{!}{%
		$
		\begin{array}{ccccccc}
			\cD_{1,\pm 0} \times \cD_{2,\pm 1} & \to & [\cD_{4,\pm \{1,1,1,3\}}\times \cD_{2,\pm 0}]_\even & \xrightarrow{\text{halve}} & \cD_{2,\pm \{0,1,1,1\}} \times \cD_{1,\pm 0} & \to & (\cD_{1,\pm 0} \times \cD_{2,\pm 1})^2\\[3pt]
			(\pi_\pm^{(1)}, \pi_\pm^{(2)}) & \mapsto & \big((\mu_\pm^{(1)}, \mu_\pm^{(2)}, \mu_\pm^{(3)}, \mu_\pm^{(4)}), \xi_\pm\big) & \xmapsto{\text{\phantom{halve}}} & \big((\omega_\pm^{(1)}, \omega_\pm^{(2)}, \sigma_\pm^{(2)}, \rho_\pm^{(2)}),\rho_\pm^{(1)}\big) & \mapsto & \big((\sigma_\pm^{(1)}, \sigma_\pm^{(2)}), (\rho_\pm^{(1)}, \rho_\pm^{(2)})\big)
		\end{array}
		$
	}
\end{equation*}
with the desired relation
\begin{align*}
	\len(\pi_\pm^{(2)}) = \len(\sigma_\pm^{(2)}) + \len(\rho_\pm^{(2)}).
\end{align*}

\section{Conclusion}\label{sec:conclusion}

We raise two questions for future study.

First, let $\dis(\lambda)$ denote the number of different part sizes in a partition $\lambda$ with designated summands. It was shown in \cite[Theorem~1.3]{FuSel} that
\begin{align}\label{eq:2-conv-PDO-dis-comb}
	\UU_2\left(\sum_{\lambda\in \cPDO} x^{\dis(\lambda)} q^{|\lambda|}\right) = \left(\sum_{\lambda\in \cPDO} x^{\dis(\lambda)} q^{|\lambda|}\right)^2.
\end{align}
The parametrization \eqref{eq:para-2-conv} is equivalent to this $2$-convolutive relation through the substitution of variables
\begin{align*}
	X \mapsto -\frac{2-x+\sqrt{x^2-4x}}{2}.
\end{align*}
However, since this substitution is not linear, there is no direct way to translate the bijection related to the set $\cD_{1,\pm 0} \times \cD_{2,\pm 1}$, as stated in Theorem~\ref{th:para-bij}, into the setting of the $\cPDO$ set. In other words, a bijective proof of \eqref{eq:2-conv-PDO-dis-comb} remains open.

Second, we have seen that the $2$-convolutivity of the eta-products (II.1)--(II.4) in Table~\ref{tab:2-conv} and their parametrization \eqref{eq:para-2-conv} is proven combinatorially with recourse to a bijection for the same master theta constant identity. In \cite[Theorem~1.2]{CEFSII}, a second $2$-convolutive parametrization encoding the eta-products (II.2), (II.4), (II.5), and (II.6) was established. Furthermore, in \cite{CEFSI}, the concept of $2$-convolutivity was generalized to generic $m$-convolutivity. Following this, a $3$-convolutive parametrization that encodes fourteen eta-products was proven in \cite[Theorem~1.4]{CEFSII}. It is then intriguing to know if there are other theta constant identities such that bijections for them would produce comprehensive treatments of the two sets of objects.

\subsection*{Acknowledgements}

Shane Chern was supported by the FWF Austrian Science Fund (grant no.~10.55776/F1002). Shishuo Fu was supported by the National Natural Science Foundation of China (grants no.~12171059 and 12671383) and the Fundamental Research Funds for the Central Universities (grant no.~2025CDJ-IAISYB-008). 

\bibliographystyle{amsplain}

\appendix

\section{Example for the bijection related to the $\cPDO$ set}\label{app:example}

The purpose of this appendix is to present a concrete example to explain the bijective proof of the $2$-convolutivity of the $\PDO$ function shown in Theorem~\ref{thm:PDOb}. For simplicity, we use the \emph{multiplicity notation} for partitions,
\begin{align*}
	1^{m_1} 2^{m_2} 3^{m_3} \cdots,
\end{align*}
where $m_i$ counts the parts of size $i$, and we further extend this notation to partitions with designated summands,
\begin{align*}
	1_{p_1}^{m_1} 2_{p_2}^{m_2} 3_{p_3}^{m_3} \cdots,
\end{align*}
where $p_i$ is the position of the designated $i$. If $i$ is not a part so that $m_i = 0$, then we do not record the corresponding term in the multiplicity notation.

\begin{example}\label{ex:main}
	The partition
	\begin{align*}
		1_{33}^{34}\; 3_{12}^{14}\; 5_1^3\; 7_2^2\; 11_1^1\; 13_1^1\; 15_1^1
	\end{align*}
	in $[\cPDO]_\even$ is mapped to
	\begin{align*}
		(1_{1}^{3}\; 3_2^6\; 5_2^2\; 7_1^1,\;\; 1_{11}^{17}\; 3_2^4\; 5_1^1)
	\end{align*}
	in $\cPDO^2$ under the bijection in Theorem~\ref{thm:PDOb}.
\end{example}

In doing so, each of Lemma~\ref{le:trival}, Theorem~\ref{thm:SunKim}, Corollary~\ref{coro:Kim}, Theorem~\ref{th:ALL-bij}, Corollary~\ref{coro:PDO-bij-step1}, and Theorem~\ref{th:PDO-bij-step2} will be illustrated, and the above example for the desired bijection in Theorem~\ref{thm:PDOb} follows by gluing these pieces together. We remark that on several occasions, different colors are assigned to components of partition tuples to keep track of their rearrangement.

\subsection{Example for Lemma~\ref{le:trival}}\label{ex:trivial}

We work on the case where $m=6$ and $(b_1,b_2,b_3,b_4,b_5,b_6)=(0,1,1,2,2,3)$. Choose the sextuple
{\footnotesize\begin{align*}
	\big({\color{brown} (12, 12)},\,{\color{teal} (2, 10)},\,{\color{violet} (14, \varnothing)},\,{\color{pink} (\varnothing, \varnothing)},\,{\color{magenta} (16+4, 8)},\,{\color{olive} (\varnothing, 6)}\big)
\end{align*}}%
in $\cD_{12,\pm \{0,2,2,4,4,6\}}$. Here each component in our partition tuple has been colored differently --- this is necessary because we have two copies of $\cD_{12,\pm 2}$ and $\cD_{12,\pm 4}$ to distinguish, respectively. Now halving the parts gives
{\footnotesize\begin{align*}
	\big({\color{brown} (6, 6)},\,{\color{teal} (1, 5)},\,{\color{violet} (7, \varnothing)},\,{\color{pink} (\varnothing, \varnothing)},\,{\color{magenta} (8+2, 4)},\,{\color{olive} (\varnothing, 3)}\big)
\end{align*}}%
in $\cD_{6,\pm \{0,1,1,2,2,3\}}$. For $b\in \{0,1,1,2,2,3\}$, we map $\cD_{6,\pm b}$ to $\cD_{12,\pm b}\times \cD_{12,\pm (6-b)}$, or more precisely, ${\color{red} \cD_{6, b}} \times {\color{blue} \cD_{6, 6-b}}$ to $({\color{red} \cD_{12, b}} \times {\color{blue} \cD_{12, 12-b}}) \times ({\color{blue} \cD_{12, 6-b}} \times {\color{red} \cD_{12, 6+b}})$. Thus, the previous sextuple is mapped to
{\footnotesize\begin{align*}
	\big({\color{brown} (\varnothing,\varnothing)},\,{\color{teal} (1,\varnothing)},\,{\color{violet} (\varnothing,\varnothing)},\,{\color{pink} (\varnothing,\varnothing)},\,{\color{magenta} (2,\varnothing)},\,{\color{olive} (\varnothing,\varnothing)},\,{\color{olive} (3,\varnothing)},\,{\color{magenta} (4,8)},\,{\color{pink} (\varnothing,\varnothing)},\,{\color{violet} (\varnothing,7)},\,{\color{teal} (5,\varnothing)},\,{\color{brown} (6,6)}\big)
\end{align*}}%
in $\cD_{12,\pm \{0,1,1,2,2,3,3,4,4,5,5,6\}}$.

\subsection{Example for Theorem~\ref{thm:SunKim} and Corollary~\ref{coro:Kim}}\label{ex:Kim}

We work on the case where $m=6$ and $(a_1,a_2,a_3) = (1,3,3)$ in Theorem~\ref{thm:SunKim} and Corollary~\ref{coro:Kim}. 

\subsubsection{Example for Theorem~\ref{thm:SunKim}}\label{ex:Kim-original}

Choose the quadruple
\begin{align*}
	\big((13+1,11),\, (15+3,\varnothing),\, (3,9),\, (\varnothing,5)\big)
\end{align*}
in $[\cD_{12,\pm \{1,3,3,7\}}]_\even$. Applying the quintipartite procedure in Section~\ref{sec:Kim-bij-2} gives
{\footnotesize\begin{align}
	\big((13+1,11),\, (15+3,\varnothing),\, (3,9),\, (\varnothing,5)\big) &\xmapsto{\text{Step (1)}} \big((1+0,1),\, (1+0,\varnothing),\, (0,1),\, (\varnothing,1)\big) \nonumber \\
	&\xmapsto{\text{Step (2)}} \big((2,\varnothing,1,\varnothing),\, (1,2,0,-1)\big) \nonumber \\
	&\xmapsto{\text{Step (3)}} \big((2,\varnothing,1,\varnothing),\, (2,1,-1,0)\big) \nonumber \\
	&\xmapsto{\text{Step (4)}} \big((2+1,\varnothing),\, (0,\varnothing),\, (\varnothing,2),\, (\varnothing,\varnothing)\big) \nonumber \\
	&\xmapsto{\text{Step (5)}} \big((24+12,\varnothing),\, (4,\varnothing),\, (\varnothing,20),\, (\varnothing,\varnothing)\big), \label{eq:quad}
\end{align}}%
which is in the desired set $[\cD_{12,\pm \{0,4,4,6\}}]_\even$.

\subsubsection{Example for Corollary~\ref{coro:Kim}}\label{ex:Kim-coro}

Choose the quadruple
\begin{align*}
	\big((13+1,11),\, (15+3,\varnothing),\, (3,9),\, (5,\varnothing)\big)
\end{align*}
in $[\cD_{12,\pm \{1,3,3,5\}}]_\even$. Interchanging the two components in the last pair produces a quadruple in $[\cD_{12,\pm \{1,3,3,7\}}]_\even$, which is exactly the input in the previous example in Section~\ref{ex:Kim-original}. Halving the resulting quadruple \ref{eq:quad} 
then gives
\begin{align*}
	\big({\color{brown} (12+6,\varnothing)},\, {\color{teal} (2,\varnothing)},\, {\color{violet} (\varnothing,10)},\, {\color{pink} (\varnothing,\varnothing)}\big)
\end{align*}
in $\cD_{6,\pm \{0,2,2,3\}}$, where we have colored the components as before. For each $a\in \{0,2,2,3\}$, we again map $\cD_{6,\pm a}$ to $\cD_{12,\pm a}\times \cD_{12,\pm (6-a)}$. Then the above quadruple goes to
\begin{align*}
	\big({\color{brown} (12,\varnothing)},\, {\color{teal} (2,\varnothing)},\, {\color{violet} (\varnothing,10)},\, {\color{pink} (\varnothing,\varnothing)},\, {\color{pink} (\varnothing,\varnothing)},\, {\color{violet} (\varnothing,\varnothing)},\, {\color{teal} (\varnothing,\varnothing)},\, {\color{brown} (\varnothing,6)}\big)
\end{align*}
in $\cD_{12,\pm \{0,2,2,3,3,4,4,6\}}$.

\subsection{Example for Theorem~\ref{th:ALL-bij} and Corollary~\ref{coro:PDO-bij-step1}}\label{ex:ALL-bij}

Here we give three examples to illustrate Theorem~\ref{th:ALL-bij} and Corollary~\ref{coro:PDO-bij-step1}, with the first one explained in full detail. In particular, the three partitions in $\cPDO$ are exactly those in Example~\ref{ex:main}.

First, choose the partition
\begin{align*}
	1_{33}^{34}\; 3_{12}^{14}\; 5_1^3\; 7_2^2\; 11_1^1\; 13_1^1\; 15_1^1
\end{align*}
in $[\cPDO]_\even$. In light of Theorem~\ref{th:ALL-bij},
\begin{align*}
	\begin{array}{rclp{2.5em}rcl}
		1_{33}^{34} &\xrightarrow[p_1=33]{m_1=34,\; m_1'=18}& (1^3, 2^{14}, 3^1), & & 3_{12}^{14} &\xrightarrow[p_3=12]{m_3=14,\; m_3'=8}& (3^5, 6^3, 9^1),\\
		5_1^3 &\xrightarrow[p_5=1]{m_5=3,\; m_5'=2}& (5^3, \varnothing, \varnothing), & &
		7_2^2 &\xrightarrow[p_7=2]{m_7=2,\; m_7'=2}& (\varnothing, 14^1, \varnothing),\\
		11_1^1 &\xrightarrow[p_{11}=1]{m_{11}=1,\; m_{11}'=1}& (11^1, \varnothing, \varnothing), & &
		13_1^1 &\xrightarrow[p_{13}=1]{m_{13}=1,\; m_{13}'=1}& (13^1, \varnothing, \varnothing),\\
		15_1^1 &\xrightarrow[p_{15}=1]{m_{15}=1,\; m_{15}'=1}& (15^1, \varnothing, \varnothing).
	\end{array}
\end{align*}
Thus, our partition in $\cPDO$ is mapped to the triple
\begin{align*}
	(1^3\, 3^5\, 5^3\, 11^1\, 13^1\, 15^1,\;\; 2^{14}\, 6^3\, 14^1,\;\; 3^1\, 9^1)
\end{align*}
in $\cP_{\mathbb{O}}\times 2\cP_{\mathbb{O}}\times 3\cD_{\mathbb{O}}$. Using Glaisher's bijection, we see that $1^3\, 3^5\, 5^3\, 11^1\, 13^1\, 15^1$ is mapped to $15+13+12+11+10+5+3+2+1$ and that $1^{14}\, 3^3\, 7^1$ is mapped to $8+7+6+4+3+2$. We conclude from Corollary~\ref{coro:PDO-bij-step1} that the above triple is mapped to
\begin{align*}
	(15+13+12+11+10+5+3+2+1,\, 16+14+12+8+6+4,\, 9+3),
\end{align*}
which is in $\cD\times \cD_{2,0}\times \cD_{6,3}$.

Next, we have that the partition
\begin{align*}
	1_{1}^{3}\; 3_2^6\; 5_2^2\; 7_1^1
\end{align*}
in $\cPDO$ is mapped to
\begin{align*}
	(1^{3}\, 3^4\, 7^1,\;\; 6^1\, 10^1,\;\; \varnothing)
\end{align*}
in $\cP_{\mathbb{O}}\times 2\cP_{\mathbb{O}}\times 3\cD_{\mathbb{O}}$ by Theorem~\ref{th:ALL-bij}, and hence to
\begin{align*}
	(12+7+2+1,\, 10+6,\, \varnothing)
\end{align*}
in $\cD\times \cD_{2,0}\times \cD_{6,3}$ by Corollary~\ref{coro:PDO-bij-step1}.

Finally, we have that the partition
\begin{align*}
	1_{11}^{17}\; 3_2^4\; 5_1^1
\end{align*}
in $\cPDO$ is mapped to
\begin{align*}
	(1^{12}\, 3^2\, 5^1,\;\; 2^1\, 6^1,\;\; 3^1)
\end{align*}
in $\cP_{\mathbb{O}}\times 2\cP_{\mathbb{O}}\times 3\cD_{\mathbb{O}}$ by Theorem~\ref{th:ALL-bij}, and hence to
\begin{align*}
	(8+6+5+4,\, 6+2,\, 3)
\end{align*}
in $\cD\times \cD_{2,0}\times \cD_{6,3}$ by Corollary~\ref{coro:PDO-bij-step1}.

\subsection{Example for Theorem~\ref{th:PDO-bij-step2}}\label{ex:PDO-bij-step2}

Choose the triple
{\footnotesize\begin{align*}
	({\color{blue} 15+13+12+11+10+5+3+2+1},\, {\color{olive} 16+14+12+8+6+4},\, {\color{violet} 9+3})
\end{align*}}%
in $[\cD\times \cD_{2,0}\times \cD_{6,3}]_\even$, which is related to the first example in Section~\ref{ex:ALL-bij}. We separate the parts according to their residues modulo $12$ to get
{\footnotesize\begin{align*}
	\big(({\color{blue} 12}, {\color{olive} 12}),\,({\color{blue} 13+1}, {\color{blue} 11}),\,({\color{blue} 2}, {\color{blue} 10}),\,({\color{olive} 14}, {\color{olive} \varnothing}),\,({\color{blue} 15+3}, {\color{blue} \varnothing}),\,({\color{violet} 3}, {\color{violet} 9}),\,({\color{blue} \varnothing}, {\color{blue} \varnothing}),\,({\color{olive} 16+4}, {\color{olive} 8}),\,({\color{blue} 5}, {\color{blue} \varnothing}),\,({\color{blue} \varnothing}, {\color{olive} 6})\big)
\end{align*}}%
in $\cD_{12,\pm \{0,1,2,2,3,3,4,4,5,6\}}$, and further a pair in $[\cD_{12,\pm \{1,3,3,5\}}]_\even\times \cD_{12,\pm \{0,2,2,4,4,6\}}$, with the component from $[\cD_{12,\pm \{1,3,3,5\}}]_\even$ being
{\footnotesize\begin{align*}
	\big(({\color{blue} 13+1}, {\color{blue} 11}),\, ({\color{blue} 15+3,\varnothing}),\, ({\color{violet} 3}, {\color{violet} 9}),\, ({\color{blue} 5}, {\color{blue} \varnothing})\big),
\end{align*}}%
and the component from $\cD_{12,\pm \{0,2,2,4,4,6\}}$ being
{\footnotesize\begin{align*}
	\big(({\color{blue} 12}, {\color{olive} 12}),\, ({\color{blue} 2}, {\color{blue} 10}),\, ({\color{olive} 14}, {\color{olive} \varnothing}),\, ({\color{blue} \varnothing}, {\color{blue} \varnothing}),\, ({\color{olive} 16+4}, {\color{olive} 8}),\, ({\color{blue} \varnothing}, {\color{olive} 6})\big).
\end{align*}}%
According to the example in Section~\ref{ex:Kim-coro}, the weight-halving bijection in Corollary~\ref{coro:Kim} maps the previous quadruple to
{\footnotesize\begin{align*}
	\big({\color{teal} (12,\varnothing)},\,{\color{teal} (2,\varnothing)},\,{\color{teal} (\varnothing,10)},\,{\color{teal} (\varnothing,\varnothing)},\,{\color{teal} (\varnothing,\varnothing),\,{\color{teal} (\varnothing,\varnothing)},\,{\color{teal} (\varnothing,\varnothing)},\,{\color{teal} (\varnothing,6)}}\big)
\end{align*}}%
in $\cD_{12,\pm \{0,2,2,3,3,4,4,6\}}$; according to the example in Section~\ref{ex:trivial}, the weight-halving bijection in Lemma~\ref{le:trival} maps the previous sextuple to
{\footnotesize\begin{align*}
	\big({\color{BlueViolet} (\varnothing,\varnothing)},\,{\color{orange} (1,\varnothing)},\,{\color{BlueViolet} (\varnothing,\varnothing)},\,{\color{BlueViolet} (\varnothing,\varnothing)},\,{\color{BlueViolet} (2,\varnothing)},\,{\color{BlueViolet} (\varnothing,\varnothing)},\,{\color{BlueViolet} (3,\varnothing)},\,{\color{BlueViolet} (4,8)},\,{\color{BlueViolet} (\varnothing,\varnothing)},\,{\color{orange} (\varnothing,7)},\,{\color{BlueViolet} (5,\varnothing)},\,{\color{BlueViolet} (6,6)}\big)
\end{align*}}%
in $\cD_{12,\pm \{0,1,1,2,2,3,3,4,4,5,5,6\}}$. Here we have introduced a second coloring system for the rearrangement toward an ordered pair in $(\cD_{12,\pm \{0,1,2,2,3,3,4,4,5,6\}})^2$, with the first component being
{\footnotesize\begin{align*}
	\big({\color{teal} (12,\varnothing)},\,{\color{orange} (1,\varnothing)},\,{\color{teal} (2,\varnothing)},\,{\color{teal} (\varnothing,10)},\,{\color{teal} (\varnothing,\varnothing)},\,{\color{teal} (\varnothing,\varnothing)},\,{\color{teal} (\varnothing,\varnothing)},\,{\color{teal} (\varnothing,\varnothing)},\,{\color{orange} (\varnothing,7)},\,{\color{teal} (\varnothing,6)}\big)
\end{align*}}%
and the second component being
{\footnotesize\begin{align*}
	\big({\color{BlueViolet} (\varnothing,\varnothing)},\,{\color{BlueViolet} (\varnothing,\varnothing)},\,{\color{BlueViolet} (\varnothing,\varnothing)},\,{\color{BlueViolet} (2,\varnothing)},\,{\color{BlueViolet} (\varnothing,\varnothing)},\,{\color{BlueViolet} (3,\varnothing)},\,{\color{BlueViolet} (4,8)},\,{\color{BlueViolet} (\varnothing,\varnothing)},\,{\color{BlueViolet} (5,\varnothing)},\,{\color{BlueViolet} (6,6)}\big).
\end{align*}}%
Finally, we return to the original coloring system and pull $\cD_{12,\pm \{0,1,2,2,3,3,4,4,5,6\}}$ back to $\cD\times \cD_{2,0}\times \cD_{6,3}$. Then
{\footnotesize\begin{align*}
		\big(({\color{blue} 12}, {\color{olive} \varnothing}),\,({\color{blue} 1}, {\color{blue} \varnothing}),\,({\color{blue} 2}, {\color{blue} \varnothing}),\,({\color{olive} \varnothing}, {\color{olive} 10}),\,({\color{blue} \varnothing}, {\color{blue} \varnothing}),\,({\color{violet} \varnothing}, {\color{violet} \varnothing}),\,({\color{blue} \varnothing}, {\color{blue} \varnothing}),\,({\color{olive} \varnothing}, {\color{olive} \varnothing}),\,({\color{blue} \varnothing}, {\color{blue} 7}),\,({\color{blue} \varnothing}, {\color{olive} 6})\big)
\end{align*}}%
goes to the first triple
{\footnotesize\begin{align*}
	({\color{blue} 12+7+2+1},\, {\color{olive} 10+6},\, {\color{violet} \varnothing})
\end{align*}}%
in $\cD\times \cD_{2,0}\times \cD_{6,3}$, and
{\footnotesize\begin{align*}
	\big(({\color{blue} \varnothing}, {\color{olive} \varnothing}),\,({\color{blue} \varnothing}, {\color{blue} \varnothing}),\,({\color{blue} \varnothing}, {\color{blue} \varnothing}),\,({\color{olive} 2}, {\color{olive} \varnothing}),\,({\color{blue} \varnothing}, {\color{blue} \varnothing}),\,({\color{violet} 3}, {\color{violet} \varnothing}),\,({\color{blue} 4}, {\color{blue} 8}),\,({\color{olive} \varnothing}, {\color{olive} \varnothing}),\,({\color{blue} 5}, {\color{blue} \varnothing}),\,({\color{blue} 6}, {\color{olive} 6})\big)
\end{align*}}%
goes to the second triple
{\footnotesize\begin{align*}
	({\color{blue} 8+6+5+4},\, {\color{olive} 6+2},\, {\color{violet} 3})
\end{align*}}%
in $\cD\times \cD_{2,0}\times \cD_{6,3}$. It is notable that the two outputs are related to the last two examples in Section~\ref{ex:ALL-bij}, respectively.

\end{document}